\RequirePackage{fix-cm}
\documentclass[12pts,letter]{article}       
\usepackage{graphics}
\usepackage{epstopdf}
\usepackage{amssymb,amsmath}
\usepackage{graphicx}
\begin{document}

\title{Real-Time Second-Order Optimal Guidance Strategies for Optimizing Aircraft Performance in Stochastic Wind Conditions
}


\author{Kamran Turkoglu \\
(\emph{submitted and under review in Aerospace Science and Technologies})
}





\maketitle

\begin{abstract}
This study presents a real-time guidance strategy for an unmanned aerial vehicles (UAVs) that can be used to enhance their flight endurance by utilizing {\sl insitu} measurements of wind speeds and wind gradients. In these strategies, periodic adjustments are made in the airspeed and/or heading angle command, in level flights, for the UAV to minimize a projected power requirement. In this study, UAV dynamics are described by a three-dimensional dynamic point-mass model.
A stochastic wind field model has been used to analyze the effect of the wind in the process. 
Onboard closed-loop trajectory tracking logics that follow airspeed vector commands are modeled using the method of feedback linearization. To evaluate the benefits of these strategies in enhancing UAV flight endurance, a reference strategy is introduced in which the UAV would follow the optimal airspeed command in a steady level flight under zero wind conditions. A performance measure is defined as the average power consumption with respect to no wind case. Different scenarios have been evaluated both over a specified time interval and over different initial heading angles of the UAV. A relative benefit criterion is then defined as the percentage improvement in the performance measure of a proposed strategy over that of the reference strategy. Extensive numerical simulations are conducted to show efficiency and applicability of the proposed algorithms. Results demonstrate possible power savings of the proposed real-time guidance strategies in level flights, by utilization of wind energy.
\end{abstract}

\section{Introduction}\label{intro}
When nature is observed, it is possible to see that birds, are capable of taking advantage of wind currents not only to minimize their energy consumption, but also to maximize their endurance. One important aspect of this is, they do not hold any information and the ability to estimate the weather conditions on the path they fly (and/or going to fly through). All decisions are solely based on with respect to local and instantaneous wind conditions. In this way, they optimize flight trajectory based on the local and instantaneous decisions. This is, a simple but, yet, inspiring mechanism to learn and incorporate into flight dynamics through the mechanics of flight.

Ideally, if the regional wind information is completely known in advance (with the help of a pre-determined (forecasted) weather/wind maps over the flight region), optimal flight/trajectory planning can be used to determine flight paths that minimize the total power consumption over a specified time interval, subject to various constraints. However, the main challenge for such a \emph{pre-determined map} approach is that (due to the highly complex, stochastic, coupled and  nonlinear nature of the atmosphere) weather forecasting related prediction errors also propagate into the optimization routine.

Therefore, instead of making decisions based on pre-determined weather maps, with this study, we propose real-time guidance strategies that will make \emph{local}, \emph{in-situ} decisions using \emph{available on-board} instruments to benefit from the existing \emph{local} wind conditions and minimize power consumption during the flight. The foundation of this concept had been briefly outlined in Turkoglu's fundamental work \cite{Turkoglu2009,TurkogluPhd2012}.

There are pioneering works in the area of UAV flights utilizing wind energies. The developments and flight tests of practical guidance strategies for detecting and utilizing thermals \cite{Allen2005,Allen2007}
have illustrated the feasibility of these concepts.
Patel \cite{Patel2006} studied the effect of wind in determining optimal flight control conditions under the influence of atmospheric turbulence. Langelaan \cite{Langelaan2008} studied how to exploit energy from high frequency gusts in the vertical plane for UAVs.
In addition, Langelaan \cite{Langelaan2009} presented a method for minimum energy path planning in complex wind fields using a predetermined energy map. Sukkarieh et al. \cite{Sukkarieh2009a,Sukkarieh2009b,SukkariehLawrence2011} developed a framework for an energy-based path planning that utilizes local wind estimations for dynamic soaring using the measurements and predictions from the wind patterns. Rysdyk \cite{Rysdyk2007} studied the problem of course and heading changes in significant wind conditions. McNeely \cite{McNeely2007} and et al. studied the tour planning problem for UAVs under wind conditions. McGee \cite{McGee2007} presented a study of optimal path planning using a kinematic aircraft model. 

In all of these studies, it is assumed that wind information is fully known over the region of flight. But in reality, wind is a stochastic process that needs to be addressed accordingly. One suitable approach is to devise real-time strategies that will benefit from the instantaneous nature of wind dynamics, on the spot, rather than depending on big forecasted weather maps. Thus, it is reasonable to adopt the idea of executing \emph{local}, \emph{instantaneous} maneuvers based on local/on-board measurements to utilize wind energy via the information available at that specific time instant $t_0$. Compared with the dynamic optimization studies, this study presents the use and the utilization of \emph{in-situ} wind measurements alone with no regional wind information, to optimize power consumptions.

In this paper, optimal adjustments are made to the airspeed, heading angle and/or flight path angle commands to minimize a projected power consumption, based on the \emph{instantaneous}, \emph{local} wind conditions/measurements. The onboard feedback control system then tracks these modified (updated) commands, and this process is repeated periodically throughout the entire flight.

\section{Dynamic system analysis}
\label{ch:uav_model_Section}

When the nature of trajectory optimization problems is taken into account, it is possible to see that there are six main components which are vital in determining (and also analyzing) the flight trajectory: airspeed (V), heading angle ($\Psi$), flight path angle ($\gamma$) and location of the aircraft, namely $x$, $y$ and $h$. Once these values are known (and/or provided), trajectory planning becomes a relatively easy task.

In this study, for the purpose of developing optimal guidance strategies, UAV flights are represented using a 3D dynamic-point-mass model, and the detailed structure is provided in further detail in the following sections.

\subsection{Normalized equations of motion}
 In order to increase numerical efficiency and to reduce computational complexity, 3-D dynamic point mass equations of motion are normalized by specifying a characteristic air-speed $V_n$ and mass $m$. In this paper, characteristic \emph{normalization speed-$V_n$} is selected to be the maximum speed of the aircraft (i.e. $V_n = V_{max}$).

Following to some algebraic manipulations \cite{TurkogluPhd2012}, it is possible to obtain normalized equations of motion as
\begin{equation}\label{VdotNormalized}
\begin{split}
&\bar{V}' = {\bar{P} \over \bar{V}} - \bar{\rho} \bar{V}^2 (C_{D_0} + K C^2_L) - \sin\gamma - \bar W_V'(\bar{V}, \Psi, \gamma, \bar{x}, \bar{y}, \bar{h})
\end{split}
\end{equation}
\begin{equation}\label{PsidotNormalized}
\begin{split}
& \Psi' = {{\bar{\rho} \bar V C_L}\over\cos\gamma}\sin\mu - {1\over \bar V \cos\gamma} \bar W_\Psi'(\bar{V}, \Psi, \gamma, \bar{x}, \bar{y}, \bar{h})
\end{split}
\end{equation}
\begin{equation}\label{GammadotNormalized}
\begin{split}
& \gamma'  = \bar{\rho} \bar{V} C_L\cos\mu- {\cos\gamma\over \bar{V}}+ {1 \over{\bar{V}}} \bar W_\gamma'(\bar{V}, \Psi, \gamma, \bar{x}, \bar{y}, \bar{h})
\end{split}
\end{equation}
\begin{equation}\label{XdotNormalized}
\bar{x}'  = \bar{V} \cos\gamma \sin\Psi + \bar{W}_x(\bar{x}, \bar{y}, \bar{h})
\end{equation}
\begin{equation}\label{YdotNormalized}
\bar{y}'  = \bar{V} \cos\gamma \cos\Psi + \bar{W}_y(\bar{x}, \bar{y}, \bar{h})
\end{equation}
\begin{equation}\label{HdotNormalized}
\bar{h}'  = \bar{V} \sin\gamma + \bar{W}_h(\bar{x}, \bar{y}, \bar{h})
\end{equation}
and
\begin{equation}\label{eq:WvdotNormalized_new_last}
\bar{W}_V' = \bar{W}'_x \cos\gamma \sin\Psi + \bar{W}'_y \cos\gamma\cos\Psi + \bar{W}'_h \sin\gamma \\
\end{equation}
\begin{equation}\label{eq:WpsidotNormalized_new_last}
\bar{W}_{\Psi}' = \bar{W}_x' \cos(\Psi) - \bar{W}_y' \sin(\Psi) \\
\end{equation}
\begin{equation}\label{eq:WgammadotNormalized_new_last}
\bar{W}_{\gamma}' =  \bar{W}_x' \sin(\gamma) \sin(\Psi) + \bar{W}_y' \sin(\gamma) \cos(\Psi) - \bar{W}_h' \cos(\gamma) \\
\end{equation}
where the functional dependencies of the wind terms are shown in parenthesis, for convenience.


Since this is a constrained optimization problem, imposed constraints on states and controls are also expressed using normalized values, and are presented in the following Sections.

\section{Problem statement}\label{ch:problem_statement_Section}
Guidance strategies, in general, cold be grouped into three basic categories:  {\bf action strategy}, {\bf velocity strategy}, and {\bf trajectory strategy}. 

In this study, {\bf action strategies} refer to the \emph{direct specifications} of some (or all) of the control variables, namely $z_U$ = $(P, C_L, \mu)$ (i.e. power, lift coefficient, bank angle) over a certain period of time. Here, they represent open-loop control schemes. In comparison, {\bf velocity strategies} specify some or all of the desired velocity components $ z_V = (V_c, \Psi_c, \gamma_c)$ (i.e. airspeed, heading angle, flight path angle) over a certain time interval ($\Delta t$) as flight commands, where these commands are then followed via closed-loop tracking. Finally as its name suggests, {\bf trajectory strategies} specify a flight trajectory of desired positions as functions of time over a certain time interval: $x_c(t), y_c(t), h_c(t)$ where  $t \in [t_0, t_f]$. These different categories of guidance strategies can be used to harvest wind energies of different frequencies or types, and/or be used in combination.

A {\bf complete} guidance strategy consistently needs to specify three flight commands at a given time, respectively along the longitudinal, lateral, and vertical direction. The guidance command along each direction can be from any of the three groups.

\begin{figure}[htbp]
\centering
\includegraphics[scale=0.8]{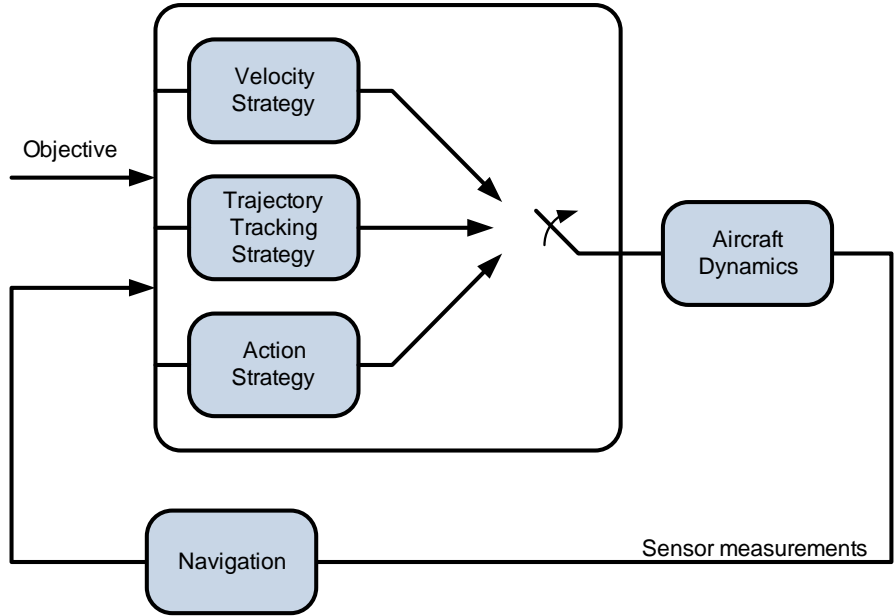}
\caption{High-level architecture of guidance strategies.}
\label{fig:highlevel_strat}
\end{figure}

Here, a high-level picture of a complete guidance strategy is provided in Fig-\ref{fig:highlevel_strat}.

In this study, it is desired to study \emph{{\bf velocity guidance strategies}}. The combination of the listed strategies above is stil an ongoing research effort, and results will be reported in another studu. 


\subsection{Problem formulation of real-time guidance strategies:}

Here, with the utilization of wind energy and velocity guidance strategy, it is desired to determine optimal adjustments in airspeed ($\bar{V}$), heading angle ($\Psi$) and/or flight path angle($\gamma$) that will minimize the power consumption projected sometime into the future. In other words, we aim to minimize the cost function expressed at the terminal cost based on local, in-situ wind measurements. Mathematically speaking: $\min_{z_{V}} \; \bar{I} = \bar{P}(\bar{t}_0+\Delta \bar{t})$, which is subject to all applicable constraints. As a result, a projected power consumption expressed at the terminal state (i.e. at $(\bar{t}_0+\Delta \bar{t})$) is used to be minimized instead of the current power at initial state (i.e. $(\bar{t}_0)$).

Feasible UAV flights must satisfy constraints due to UAV performance and operational limits. These constraints can typically be expressed as bounds on trajectory states and controls. For the current problem, these constraints are 
$$V_{\min} \leq V \leq V_{\max}$$
$$P_{\min}  \leq P \leq P_{\max}$$
$$C_{L_{\min}} \leq C_L \leq C_{L_{\max}}$$
$$|\mu| \leq \mu_{\max}$$

Following to the optimal solutions that are attained with respect to given constraints, the UAV will be directed to track optimal state commands 
$$\bar{V}^* = \bar{V}_0 + \Delta \bar{V}_c$$
$$\Psi^* = \Psi_0 + \Delta \Psi_c$$
$$\gamma^* = \gamma_0 + \Delta \gamma_c$$
For \emph{velocity guidance strategies}, all corrections ($\Delta \bar{V}$, $\Delta \Psi$ and/or $\Delta \gamma$) are the optimal adjustments (increments or decrements) to be determined with respect to \emph{local}, \emph{instantaneous} and \emph{on-board} wind measurements. 

Once optimal adjustments in airspeed ($\Delta \bar{V}$), heading angle ( $\Delta \Psi$) and in flight path angle ($\Delta \gamma$) are obtained, it takes some finite time for the UAV to achieve the desired changes via closed-loop tracking. This is a phenomenon that has to be taken into account and it is explained in further detail in the following Section.

\section{Solution strategies for level flights} \label{ch:2D_level_flight_sol_strat}


\subsection{Level flight strategies in presence of wind}

In level flight strategies, it is assumed that the aircraft is maintaining a level flight by controlling the flight path angle as zero ($\gamma = 0$). This assumption, leads to several simplifications to the presented 3D formulation, such as $\bar{W}_{\gamma}' = 0$. It also sets vertical wind component as zero ($\bar{W}_h = 0$), by the nature of problem formulation.

With these assumptions, at steady state conditions (i.e. $\bar{V}' \approx 0 $, $\Psi' \approx 0 $ and $\gamma' \approx 0 $ ), it is possible to obtain an expression for the power at time instant $\bar{t}_0$ for level flight strategies as
\begin{equation}\label{eq:P_0_3D_initial}
\begin{split}
\bar{P} (\bar{t}_0) \approx \bar{\rho} \bar{V}^3 (C_{D_0} + KC_L^2)  +  \bar{V} \bar{W}_V'
\end{split}
\end{equation}

\begin{figure}
\centering
 \includegraphics[width=10cm]{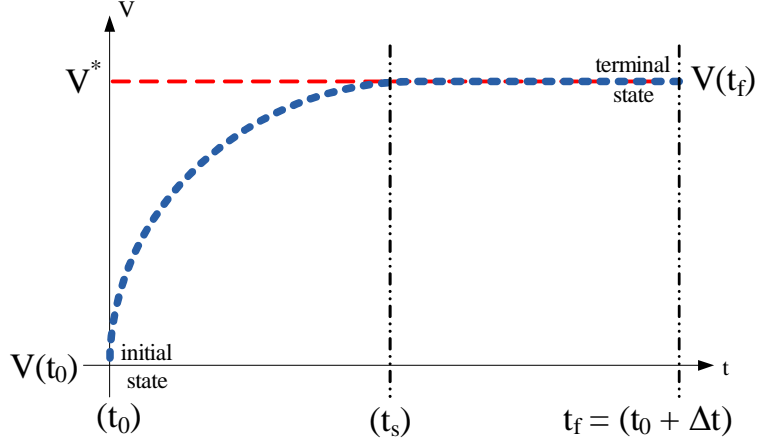}
 \caption{Trajectory tracking control logic.}\label{fig:settling_time}
\end{figure}

In 2D level flight strategies, it is assumed that by the time ($\bar{t}_0+\Delta \bar{t}$), any commanded changes in airspeed ($\bar V$) and heading angle ($\Psi$) will be achieved via closed-loop tracking, as demonstrated in Figure-\ref{fig:settling_time}. 

As a result, the vehicle will be in a steady state condition: $\bar{V}' \approx 0$, $\bar{\Psi}' \approx 0$, $\mu \approx 0$ . This leads to a control command input for lift coefficient ($C_L$) such as
\begin{equation}\label{eq:CL_level}
\begin{split}
\gamma'  \approx 0 =&~ \bar{\rho} \bar{V} C_L \cos\mu = { \cos \gamma - \bar{W}_{\gamma}' \over \bar{V}},~~C_L  = { \cos \gamma - \bar{W}_{\gamma}' \over \bar{\rho} \bar{V}^2 \cos\mu}  \\
\end{split}
\end{equation}

\noindent For level flight conditions, with the assumptions of $\gamma \approx 0$ and $\bar{W}_{\gamma}' = 0$, control command input for lift coefficient ($C_L$) in Eq.(\ref{eq:CL_level}) simplifies to,
$C_L  = { 1 / \bar{\rho} \bar{V}^2 \cos\mu}$.
Following to the derivation of control command input for lift coefficient, using Eq. (\ref{VdotNormalized}), power function is obtained as
\begin{equation}\label{eq:power}
\begin{split}
\bar{P} \approx \bar{V} \bar{V}'  + \bar{\rho} \bar{V}^3 C_{D_0} + {K  \over \bar{\rho} \cos^2\mu} {1 \over \bar{V}} + \bar{V} \bar{W}_V'  \\
\end{split}
\end{equation}
where the projected power function (at $\bar{t}_0+\Delta \bar{t}$) for level flight strategies is derived as
\begin{equation}\label{eq:minPbar3}
\begin{split}
\bar{P}(\bar{t}_0+\Delta \bar{t}) &~ \approx \bar{\rho} [\bar{V}_0 + \Delta \bar{V}_c]^3 C_{D_0} + {K  \over \bar{\rho} } {1 \over \bar{V}_0 + \Delta V_c } \\
& + [\bar{V}_0 + \Delta \bar{V}_c] \bar{W}_V'(t_0+\Delta t)
\end{split}
\end{equation}

It is assumed that by the time $\bar{t}_0+\Delta \bar{t}$, any commanded changes in airspeed and heading angle will have been mostly achieved via closed-loop tracking. Therefore, the vehicle will basically be at a steady state condition at the terminal state: $\bar{V}' \approx 0$ and $\mu \approx 0$.

Here, the $\bar{W}_V'(\bar{t}_0+\Delta \bar{t})$ term, which represents the main contribution on acceleration, needs to be expressed analytically as a function of optimal corrections and in the following section, the derivation of this term is described in further detail.

\subsection{Expression for projected wind rate in level flights}

The rate of change of any wind component can be expressed in terms of a total derivative as 
$$
\bar{W}'_{( \ )} = {\partial \bar{W}_{( \ )} \over \partial \bar{x}} \bar{x}' + {\partial \bar{W}_{( \ )} \over \partial \bar{y}} \bar{y}'
+ {\partial \bar{W}_{( \ )} \over \partial \bar{h}} \bar{h}' + {\partial \bar{W}_{( \ )} \over \partial \bar{t}}$$
where ${( \ )} = (x, y, h)$.


In level flights, with negligible flight path angle and vertical winds (i.e. $\gamma = 0$, $\bar{W}_h \approx 0$), $\bar{W}_V'$ term reduces to
\begin{equation}\label{eq:WvBar2}
\begin{split}
\bar{W}_V' = & \left( { \partial{\bar{W}_x} \over \partial{\bar{y}}} + { \partial{\bar{W}_y} \over \partial{\bar{x}}}\right) \bar{V} \sin\Psi \cos\Psi \\
& + { \partial{\bar{W}_x} \over \partial{\bar{x}}} \bar{V} \sin^2\Psi + { \partial{\bar{W}_y} \over \partial{\bar{y}}} \bar{V}  \cos^2\Psi  \\
& + \left( \bar{W}_x { \partial{\bar{W}_x} \over \partial \bar{x} } + \bar{W}_y { \partial{\bar{W}_x} \over \partial \bar{y} } + { \partial{\bar{W}_x} \over \partial \bar{t} } \right) \sin\Psi \\
& + \left( \bar{W}_x { \partial{\bar{W}_y} \over \partial \bar{x} } + \bar{W}_y { \partial{\bar{W}_y} \over \partial \bar{y} } + { \partial{\bar{W}_y} \over \partial \bar{t} } \right) \cos\Psi \\
\end{split}
\end{equation}

Since, the only available information is the \emph{local, insitu} wind measurements, it is assumed that the current wind gradients will stay constant \cite{TiwariGoshal2005} over the immediate neighborhood around the current position of the UAV, leading to an expression of $\bar{W}_V'$ at ($\bar{t}_0 + \Delta \bar{t}$), which is the projection of the wind gradient. The validity of the constant wind gradient assumption is strengthened even further with faster update rate, (i.e. smaller update time, $\Delta \bar{t}$).

With the \emph{constant wind gradient assumption} over $[\bar{t}_0, \bar{t}_0 + \Delta \bar{t}]$, we also have
\begin{equation}\label{eq:windspeeds}
\begin{split}
\bar{W}_x (\bar{t}_0+\Delta \bar{t}) & \approx \bar{W}_{x_0} + { \partial{\bar{W}_x} \over \partial \bar{x} } \Delta \bar{x} + { \partial{\bar{W}_x} \over \partial \bar{y} } \Delta \bar{y} + { \partial{\bar{W}_x} \over \partial \bar{t} } \Delta \bar{t} \\
\bar{W}_y (\bar{t}_0+\Delta \bar{t}) & \approx \bar{W}_{y_0} + { \partial{\bar{W}_y} \over \partial \bar{x} } \Delta \bar{x} + { \partial{\bar{W}_y} \over \partial \bar{y} } \Delta \bar{y} + { \partial{\bar{W}_y} \over \partial \bar{t} } \Delta \bar{t} \\
\end{split}
\end{equation}

These expressions depend on $(\Delta \bar{x}, \Delta \bar{y})$, which also have a dependency on $\Delta \bar{V}$, $\Delta \Psi$, and $\Delta \bar{t}$ through Eqs.(\ref{XdotNormalized})-(\ref{HdotNormalized}), and reciprocally on the wind components over the specified time interval. Now, in order to complete the derivation of the projected wind rate expression for level flight strategies, expressions for $(\Delta \bar{x}, \Delta \bar{y})$ will be developed in the following section.

\subsection{Expressions for position changes in level flights}

The concentration in this section will mainly be focused on developing expressions for $(\Delta \bar{x}, \Delta \bar{y})$. Their dependencies, explicitly, show on the increments of airspeed and heading angle. From Eqs.(\ref{XdotNormalized})-(\ref{YdotNormalized}), and for $\gamma=0$, it is possible to have
\begin{equation}\label{eq:DxDyTrapezoidRule}
\begin{split}
\Delta \bar{x} = \bar{x}(\bar{t}_0+\Delta \bar{t}) - \bar{x}(\bar{t}_0) &  = \int_{\bar{t}_0}^{\bar{t}_0+\Delta \bar{t}} \left(\bar{V} \sin\Psi + \bar{W}_{x} \right) d \bar{t} \\
\Delta \bar{y} = \bar{y}(\bar{t}_0+\Delta \bar{t}) - \bar{y}(\bar{t}_0) &  = \int_{\bar{t}_0}^{\bar{t}_0+\Delta \bar{t}} \left( V \cos\Psi + \bar{W}_y \right) d \bar{t} \\
\end{split}
\end{equation}
With the assumption that both airspeed ($\bar{V}$) and heading angle ($\Psi$) will have achieved their commanded values (i.e. reach steady state conditions) at the end of the specified time interval, and keeping in mind that wind speeds are obtained using Eq.(\ref{eq:windspeeds}) (through the trapezoidal rule), the numerical integration of the above equations leads to

\begin{equation}\label{eq:DxDyComplete}
\begin{split}
\Delta \bar{x} & = {1 \over Q_{level}} \left[ \left({ 2 \over \Delta \bar{t}}-{\partial{\bar{W}_{y}} \over \partial{\bar{y}}} \right) B_1 +  \left( {\partial{\bar{W}_{x}} \over \partial{\bar{y}}} \right) B_2 \right] \\
 \Delta \bar{y} & = {1 \over Q_{level}} \left[ \left({\partial{\bar{W}_{y}} \over \partial{\bar{x}}} \right) B_1 +  \left( { 2 \over \Delta \bar{t}}-{\partial{\bar{W}_{x}} \over \partial{\bar{x}}}  \right) B_2 \right]
\end{split}
\end{equation}
where
\begin{equation}\label{eq:Q_2D}
\begin{split}
Q_{level} =&~  \left( {2 \over \Delta \bar{t}} \right)^2 - {2 \over \Delta \bar{t}} \left( {\partial{\bar{W}_{x}} \over \partial{\bar{x}}}  +{\partial{\bar{W}_{y}} \over \partial{\bar{y}}} \right) \\
&+ \left( {\partial{\bar{W}_{x}} \over \partial{\bar{x}}} {\partial{\bar{W}_{y}} \over \partial{\bar{y}}} - {\partial{\bar{W}_{x}} \over \partial{\bar{y}}} {\partial{\bar{W}_{y}} \over \partial{\bar{x}}} \right)
\end{split}
\end{equation}
For a sufficiently small $\Delta t$, expression in Eq.(\ref{eq:Q_2D}) shall always be nonzero; ensuring the existence of solutions for the position change expressions.

\subsection{Guidance algorithms for level flights}

Based on the derivations above, the power consumption at terminal state (i.e. at time $\bar{t}_0 + \Delta \bar{t}$) can now be expressed as a function of the current command adjustments in airspeed and heading angle. Then, the problem of reducing future power consumptions reduces to determine $\Delta \bar{V}$ and $\Delta \Psi$ from
\begin{equation}\label{eq:StaticOptim}
\min_{\Delta \bar{V}, \Delta \Psi} \bar{P}_{t_f}  = \bar{P}(\bar{t}_0 + \Delta \bar{t}) = \bar{I} (\Delta \bar{V}, \Delta \Psi; \Delta \bar{t}, \bar{X}_0)
\end{equation}
subject to
$\Delta \bar{V}_{\max}$ and $\Delta \bar{\Psi}_{\max}$ incremental constraints, and the initial state conditions required $\bar{X}_0$.


\section{Optimization strategies}\label{ch:opt_strat_Section}

From the given problem formulation, it is possible to see that the nature of the problem is very complex, extensively coupled and highly non-linear. At this point, different algorithms may be used to solve the static optimization problem in hand. But most of the existing numerical methods heavily depend on iteration routines which are not desirable in real-time applications. They impose major drawbacks from the perspective of computation time and ``convergence rate in allowed computation time'' (which can be extremely small in some cases and applications). To avoid this to certain extent, in this Section, it is aimed to solve the static optimization problem in an \emph{analytical}, and \emph{a-single-shot} manner.

For this purpose, here, yet simple but powerful \emph{gradient method} will be presented to aid in solving the static optimization problem in hand.

\subsection{Second order optimal adjustment strategies}
With second-order gradient algorithms, the main goal is to find \emph{locally} optimal adjustments: $\Delta \bar V$, $\Delta \Psi$ and $\Delta \gamma$ using the necessary conditions for optimality. For this purpose, second order Taylor series approximation of projected power function, in the neighbourhood of initial state $\bar{X}_0$, is taken into account


\normalsize
\noindent where the analytical expressions for the \emph{local minima in the neighbourhood of in-situ measurements} are obtained through the necessary conditions of optimality,

\noindent General matrix representation of problem solution can be given as $
\left[ {T_1} \right]_{(nx1)} + \left[ {T_2} \right]_{(nxn)} \left[{\Delta}\right]_{(nx1)} = 0$,
where for this case, corresponding components are obtained as
\begin{equation}
\begin{split}
 \Delta =&~
\left[ \begin{array}{c}
\Delta \bar{V} \\ \Delta \Psi \\ \Delta \gamma
\end{array} \right], ~~~
T_1 = \left[ \begin{array}{c}
 \left( {\partial \bar I \over \partial \bar{V}^*}\right)_0 \\
 \left( {\partial \bar I \over \partial \Psi^*}\right)_0 \\
 \left( {\partial \bar I \over \partial \gamma^*}\right)_0 \end{array} \right],~~~ \\
T_2 =&~ \left[
\begin{array}{ccc}
\left( {\partial^2 \bar I \over \partial \bar{V^*}^2}\right)_0 & \left( {\partial^2 \bar I \over \partial \bar{V}^* \partial \Psi^*}\right)_0 & \left( {\partial^2 \bar I \over \partial \bar{V}^* \partial \gamma^*}\right)_0 \\
\left( {\partial^2 \bar I \over \partial \Psi^* \partial \bar{V}^*}\right)_0  & \left( {\partial^2 \bar I \over \partial {\Psi^*}^2}\right)_0 & \left( {\partial^2 \bar I \over \partial \Psi^* \partial \gamma^*}\right)_0 \\
\left( {\partial^2 \bar I \over \partial \gamma^* \partial \bar{V}^*}\right)_0 & \left( {\partial^2 \bar I \over \partial \gamma^* \partial \Psi^*}\right)_0 & \left( {\partial^2 \bar I \over \partial {\gamma^*}^2}\right)_0
\end{array}
\right]
\end{split}
\end{equation}
\noindent It is possible to rewrite solution as $\left[ {T_1} \right]_{(nx1)} - \left[ - {T_2} \right]_{(nxn)} \left[{\Delta}\right]_{(nx1)} = 0$. This result leads to a well known least squares minimization problem of
$$\inf_{x \in \emph{X}} ||~ y - Ax ~||$$ 
where the optimal solution is obtained as $x^{*} = \left( A^T A \right)^{-1}A^T y$.

Therefore, in presence of favourable wind conditions, the solution to second-order optimal adjustment strategy ($\Delta \bar{V}$, $\Delta \Psi$ and $\Delta \gamma$) for the case of 2D flight is obtained as $\left[ \Delta \right]_{nx1}^{opt} = \left[ \left( T_2^T \right)_{nxn} \left( T_2 \right)_{nxn} \right]^{-1} \left( -T_2^T \right)_{nxn} \left( T_1 \right)_{nx1}$, where $(~)^T$ and $(~)^{-1}$ defines the transpose and inverse of a matrix, respectively. For the 2D level flight case, this solution reduces to calculating only $\Delta \bar{V}$ and $\Delta \Psi$. 

This completes the derivation of optimal solution strategies and provides an \emph{analytical}, \emph{single-shot-solution} to the static optimization problem in hand, in presence of wind.
Once optimal adjustments ($\Delta \bar{V}$, $\Delta \Psi$ and $\Delta \gamma$) are obtained, existing states are updated instantaneously
$$\bar{V}_c^* = \bar{V}_0 + \Delta \bar{V}$$
$$\psi_c^* = \Psi_0 + \Delta \Psi$$
and optimal flight conditions are calculated. Using appropriate control strategies, tracking calculated optimal commands 
leads to a locally optimal flight condition that minimizes power consumption in presence of wind. Next section discusses such tracking control commands methodologies.



\section{Models of closed-loop tracking}
\label{ch:models_closed_loop_tracking_Section}

In literature, there are various (general) models of control systems to achieve tracking performance. One of them is the well known \emph{feedback linearization method}. The method of feedback linearization is a widely used and a commonly accepted application in applied nonlinear control, and is well defined in many sources, as \cite{SlotineLi1991}. In the following subsections, necessary control command inputs will be derived using \emph{the method of feedback linearization} which will make sure that such desired states ($\bar{V}_c$, $\Psi_c$, $\gamma_c$) will be achieved and tracked accordingly, to minimize power consumption in presence of wind.

It should be noted that, with respect to the assigned mission objectives, the aircraft could be guided follow \emph{velocity guidance strategies} ($\bar{V}_c$,$\Psi_c$,$\gamma_c$), \emph{trajectory strategies} ($\bar{x}_c$,$\bar{y}_c$,$\bar{h}_c$) or the specific combination of both strategies. From that perspective, it is necessary to derive such control commands, which will make sure the corresponding strategy is executed as desired. Thus, in the following sections, necessary control commands for \emph{velocity guidance strategies} will be derived in detail.

\subsection{Control commands for velocity guidance strategy ($\bar{V}_{air}$, $\Psi$, $\gamma$):}
\label{Sec:control_vel_strat}

The \emph{method of feedback linearization} simply defines a transformation from the nonlinear system in hand into an equivalent linear system through change of variables and a suitable control input. Here, we would like to acknowledge the fact that effects of plant and modelling uncertainties will have an important impact on the outcome of such methodology. Therefore, how those plant/modeling uncertainties are handled is an important matter. It is a topic which should be investigated in further detail. However, it is a topic which is out of the scope of the study presented in this paper and is left aside for future research. As a result, such uncertainty analysis is not included in this paper.

In follows, the procedure for the derivation of control inputs (\emph{power}-$\bar{P}$, lift coefficient-$C_L$ and bank angle-$\mu$) is presented in terms of normalized EoMs.

If the derivation is started with first order approximation of air-speed tracking, it is possible to define the closed-loop power law as
\begin{equation}
{\bar{P}_{(vel)} \over \bar{V}} = -\bar{K}_V (\bar{V} - \bar{V}_c)  + \bar{\rho} \bar{V}^2 \left(C_{D_0} + K C_{L_v}^2 \right) + \sin\gamma + \bar{W}_V'
\end{equation}
Similarly the closed loop control command for bank angle, ($\mu$) is derived through bank angle expressions as
\begin{equation}
\tan \mu_{(vel)} = {\bar{W}_{\Psi}' -\bar{V} cos(\gamma) \bar{K}_{\Psi} (\Psi - \Psi_c) \over \cos\gamma - \bar{W}_\gamma' -\bar{V} K_{\gamma} (\gamma - \gamma_c)}
\end{equation}
where the closed loop control logic for lift coefficient ($C_L$) is obtained as given in Eq.(\ref{eq:CL_control}).
\begin{figure*}[!t]
\normalsize
\begin{equation}\label{eq:CL_control}
\begin{split}
& C_{L_{(vel)}} = {\sqrt{\left[ \bar{W}_{\Psi}' -\bar{V} \cos\gamma \bar{K}_{\Psi} (\Psi - \Psi_c) \right]^2 + \left[ \cos\gamma - \bar{W}_{\gamma}' -\bar{V} \bar{K}_{\gamma} (\gamma - \gamma_c) \right]^2} \over \bar{\rho}\bar{V}^2}
\end{split}
\end{equation}
\hrulefill
\end{figure*}

In the given derivations above, the \emph{normalized} feedback gains $(\bar{K}_V, \bar{K}_\Psi, \bar{K}_\gamma)$ can be tailored to reflect typical closed-loop UAV control characteristics. It is also good to keep in mind that, even if control commands are derived for a general 3D case, for 2D level flight strategies they simplify further with ${\gamma}'\approx 0$.


\section{Simulation evaluation and results for level flight strategies}
\label{ch:simulation_and_results_Section}


%

\subsection{Guidance algorithm parameters}

Performance of the proposed level flight guidance strategies strongly depend on the bounds (constraints) imposed on the adjustment parameters, such as $\Delta \bar{V}_{\max},~~\Delta \bar{V}_{\min},~~\Delta \Psi_{\max},~~\Delta \Psi_{\min}$. It is crucial to come up with ranges of their appropriate values for constrained optimization problems.

Actual values used in the guidance strategies can be smaller than the given bounds, but in this simulational study, we select $\Delta \bar{V}_{\max} \leq 5$ [ft/sec], $\Delta \bar{V}_{\min}  = - \Delta \bar{V}_{\max}$, $ \Delta \Psi_{\max} \leq 30^\circ$, and $\Delta \Psi_{\min} = - \Delta \Psi_{\max}$. It is also important to note that those values are specific to an UAV and is subjected to vary from aircraft to aircraft.

\subsection{Evaluation criterion}

Because of the reason that guidance strategies provided in this study are introduced to save power in UAV flights, a basic performance measure is defined as the \emph{average power consumption} over a specified time interval,
\begin{equation}
\bar{P} = {1 \over \bar{t}_f} \int^{\bar{t}_f}_{0} \bar{T} \bar{V} d \bar{t}
\end{equation}
where $t_f$ is the time period of evaluation (i.e. flight/simulation time).

In the numerical results that are presented in the following sections, sampling time (integration step size) is selected as $50[Hz]$, (i.e. $0.02[sec]$).


Furthermore, different initial heading angles of the UAV result in different relative angles with respect to the wind field, thus affecting the overall result of power saving benefit. In order to filter out these differences caused by different initial headings of UAV flights, and to have an estimate of global benefit, the above basic performance measure is further averaged within a set of different initial heading angles over $[0, 360^o]$. Here, a generic case of $\Delta \Psi_0 = 5^{\circ}$ increments are considered. The mean of basic average power consumption over different initial heading conditions is defined as the measure of the performance
\begin{equation}\bar{P}_{\rm avg} = {1 \over N_\Psi} \sum^{N_\Psi}_{i=1} \bar{P}_{\Psi_0} \end{equation}
where each $i$ corresponds to a different initial heading angle, and $N_{\Psi}$ is the number of different initial heading angles used.

To evaluate the proposed guidance strategies, following two scenarios are considered:
\begin{itemize}
\item {\bf Scenario-0:} is the reference strategy that aims to follow the reference airspeed and a constant heading angle command set at the initial heading angle. This provides the reference average power consumption, $\bar{P}^0$, in case of no wind.
\item {\bf Scenario-1:} is the case where both airspeed and heading angle commands are adjusted periodically based on the current wind measurements. The resulting average power consumption is $\bar{P}^3$.
\end{itemize}


\subsection{Benefit criterion}

To assess the outcomes of proposed strategies, the following benefit criterion is introduced as a relative measure of potential fuel savings of the proposed guidance strategies over the reference strategy
\begin{equation} B_i = {\bar{P}^0_{\rm avg} - \bar{P}^i_{\rm avg} \over \bar{P}^0_{\rm avg} } \hskip 10pt i = 1, 2, 3, 4 \end{equation}

\subsection{Simulation parameters}

In simulations, UAV parameters similar to those of the ScanEagle UAV are used, where characteristic values, taken from \cite{Austin2010}
%
Here, the power available is assumed to be able to vary instantaneously ($\%100$ change in 1sec for propeller-driven engines). In addition, available control rate bounds are taken into account as ($\bar{P}_{max}' = {6.216}$, ($\%100$ change in 1sec),($C_{L_{max}}' = 1.865$ (${1.2 \over 4} [1/sec]$),
($\mu_{max}' = 1.085$ (10 [deg/sec]).

Furthermore, the minimum airspeed constraint is selected to be close to the stall speed, whereas the maximum speed is selected to be close to a typical cruise speed, with allowances for transient dynamics in both cases. In terms of normalized quantities,$\bar{V}_{\min} = {1.0 / \sqrt{\bar{\rho} C_{L_{\max} } } }$, $\bar{V}_{\max} = {1.5 \over \sqrt{\bar{\rho} C_{L_{\rm cr} } } }$ where $C_{L_{\rm cr}}$ is a typical cruise lift coefficient, corresponding to an angle of attack of $3^o-5^o$.

\subsection{Simulation results}
For the conducted simulations in this section, a level flight at altitude of 15,000[ft] is assumed and a total flight time of $20[min]=1200[sec]$ with a simulation step size of $0.02[sec] = 50[Hz]$ is taken into account. Normalization velocity, $V_n$, is taken as the maximum flight speed of the UAV ($V_n = 134.5 [ft/sec]$). Wind magnitude ($\bar{W}_m$) and wind direction ($\Psi_w$), obtained from the stochastic wind model formulation, taken from \cite{Turkoglu2013Ams}, is used. Wind frequency, throughout these equations is assumed to be constant and equal to $\omega_w = 0.05[rad/ft] \approx 9.4[deg/m]$. It is also assumed that periodic optimal adjustments are applied once in every $4[sec]$, which is the same as the measurement update time (i.e. $\Delta {t} = 4[sec]$).

\subsubsection{Scenario-0: Constant airspeed, constant heading angle}
In this scenario, it is desired to fly with optimal airspeed (determined for no wind case) and constant heading angle, which is the initial heading angle given in the simulation settings. With this scenario, if the mission is specified to fly with optimal airspeed and initial heading angle (say $90^o$ (with respect to true north)), UAV dynamics will execute these commands throughout the entire flight and will determine the power consumption during this specific flight routine. Obtained power consumption value is the \emph{reference} power setting, $\bar{P}_0$, and will be compared with other outcomes to measure the efficiency of proposed strategies in given flight conditions.


\subsubsection{Scenario-1: Optimal airspeed, optimal heading angle in presence of wind}
It is intuitive that combined affects of both (airspeed and heading) adjustments will help to harvest the maximum amount of benefit from the proposed guidance strategy and local wind conditions.

\begin{figure}[htbp]
\centering
\includegraphics[scale=0.45]{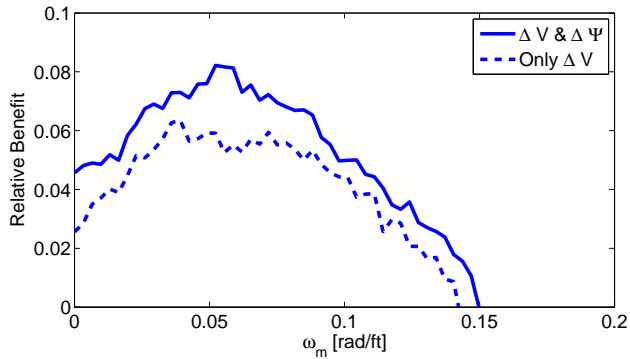}
\caption{Relative benefit versus wind frequency response of Scenario-0 and Scenario-1.}
\label{fig:RB1}
\end{figure}

For this purpose, extensive simulations have been conducted, and obtained results are given in Fig.(\ref{fig:RB1}).

Fig-\ref{fig:RB1} compares the relative benefits of the two aforementioned strategies in which periodic adjustments are made in the commands of airspeed alone (\emph{dashed line}), and both airspeed and heading (\emph{solid line}), respectively, over the reference strategy. It is possible to see from Fig-\ref{fig:RB1} that the maximum benefit is obtained for the case where adjustments both in airspeed and heading angle are applied. In this scenario, overall power savings goes up to $\%8.5$ in terms of total power consumption.

It is worhtwhile to note that for different wind and atmospheric conditions, flight path (due to heading angle adjustments) will also change. Eventually, while we seek to benefit from the wind currents and aim to extract energy from the wind, this may potentially bring the UAV to a completely undesired location, totally unrelated to the assigned flight mission. Thus, it is important to impose boundary control and restrict the flight region to complete the guidance strategies.

This concept is currently under investigation and will be reported in another study, in further details.

\section{Conclusions}
\label{ch:conclusions_chapter}

This paper presents real-time UAV guidance strategies that utilize wind energy to improve flight endurance and minimize power/fuel consumption. In these strategies, airspeed and/or heading angle commands are periodically adjusted based on the {\sl insitu} measurements of local wind components. It has been shown that using local, instantaneous wind measurements, without the knowledge of the wind field that the UAV is flying through, it is possible to benefit from the wind energy, greatly enhance the performance and increase flight endurance of the UAV. Throughout this research effort, UAV has been modeled using 3D point-mass equations. Corresponding performance and practical constraints has been introduced to mimic a realistic flight of an UAV. A stochastic wind model has been taken into account to simulate the true nature of the wind. UAV flights were formulated as a non-linear optimization problem and a cost function has been introduced to model power characteristics at terminal state as a terminal state cost function, which minimizes overall power consumption. This optimization problem has been solved as a single shot optimization, in real-time. Second-order gradient algorithms are used to find local, optimal solutions (adjustments) that will minimize the power with respect to taken local, instantaneous wind measurements. Extensive simulation results show that it is possible to obtain power savings up to $\%10$ with respect to the flight scenario, with no wind.


The proposed strategies offer improvements over the constant airspeed reference strategy in terms of average power consumptions even in a constant wind field. These benefits initially increase as the spatial frequency of the wind field gradually increases, but reaches a peak at certain frequencies and then start to decrease beyond these frequencies.

\section*{Acknowledgements}
The author is grateful to Dr.Yiyuan J. Zhao for his insightful discussions and suggestions throughout this study.

\end{document}